\documentclass[12pt]{article}
\usepackage[margin=2cm]{geometry}
\usepackage{amsfonts}
\usepackage{amssymb}
\usepackage{amsthm}
\usepackage{amsmath}
\usepackage{latexsym}
\usepackage{color}

\begin{document}
\newtheorem{lemma}{Lemma}
\newtheorem{pron}{Proposition}
\newtheorem{thm}{Theorem}
\newtheorem{Corol}{Corollary}
\newtheorem{exam}{Example}
\newtheorem{defin}{Definition}
\newtheorem{remark}{Remark}
\newcommand{\la}{\frac{1}{\lambda}}
\newcommand{\sectemul}{\arabic{section}}
\renewcommand{\theequation}{\sectemul.\arabic{equation}}
\renewcommand{\thepron}{\sectemul.\arabic{pron}}
\renewcommand{\thelemma}{\sectemul.\arabic{lemma}}
\renewcommand{\thethm}{\sectemul.\arabic{thm}}
\renewcommand{\theCorol}{\sectemul.\arabic{Corol}}
\renewcommand{\theexam}{\sectemul.\arabic{exam}}
\renewcommand{\thedefin}{\sectemul.\arabic{defin}}
\renewcommand{\theremark}{\sectemul.\arabic{remark}}
\def\REF#1{\par\hangindent\parindent\indent\llap{#1\enspace}\ignorespaces}
\def\lo{\left}
\def\ro{\right}
\def\be{\begin{equation}}
\def\ee{\end{equation}}
\def\beq{\begin{eqnarray*}}
\def\eeq{\end{eqnarray*}}
\def\bea{\begin{eqnarray}}
\def\eea{\end{eqnarray}}
\def\d{\Delta_T}
\def\r{random walk}
\def\o{\overline}
\def\P{\mathbb{P}}
\def\J{\mathcal{J}}
\def\S{\mathcal{S}}
\title{\large\bf On a transformation between distributions obeying the principle of a single big jump
\thanks{Research supported by National Science Foundation of China
(No.11071182 ).}}

\author{\small Hui Xu $^{1}$,~Michael Scheutzow $^{2}$~Yuebao Wang$^{1}$\thanks{Corresponding author.
Telephone: 86 512 67422726. Fax: 86 512 65112637. E-mail:
ybwang@suda.edu.cn}
\\
{\footnotesize\it 1. School of Mathematical Sciences, Soochow
University, Suzhou 215006, China}\\
{\footnotesize\it 2. Institut f\"ur Mathematik, Technische Universit\"at Berlin, 10623 Berlin, Germany}}
\date{}

\maketitle
{\noindent\small {\bf Abstract }}{\small Beck et al. (2013) introduced a new distribution class $\mathcal{J}$ which contains many heavy-tailed and light-tailed distributions obeying the principle of a single big jump.
Using a simple transformation which maps heavy-tailed distributions to light-tailed ones, we find some light-tailed distributions, which belong to the class $\J$ but do not belong to the convolution equivalent distribution class and which are not even weakly tail equivalent to any convolution equivalent distribution. This fact helps to understand the structure of the light-tailed distributions in the class $\J$ and leads to a negative answer to an open question raised by the above paper.
\vspace{0.2cm}

{\it Keywords:} principle of a single big jump; transformation; light-tailed distribution; convolution equivalent distribution; weak tail equivalence

\section{Introduction}
\setcounter{thm}{0}\setcounter{Corol}{0}\setcounter{lemma}{0}\setcounter{pron}{0}
\setcounter{remark}{0}\setcounter{exam}{0}\setcounter{equation}{0}

Recently, Beck et al. \cite{BBS13} introduced a new distribution class in which all distributions obey the principle of a single big jump.

Let $\{X_i,i\geq 1\}$ be a sequence of independent and identically distributed random variables with common distribution $F$.
Define the class $\J$ as the set of all distributions $F$ with unbounded support contained in $[0,\infty)$ such that for all $n\geq 2$ (or, equivalently, for $n=2$),
\begin{eqnarray}\label{101}
\lim_{K\to\infty}\liminf \P \big(X_{n,1}>x-K \mid S_n>x\big)=1,
\end{eqnarray}
where $X_{n,k}$ means the $k$-th largest random variable in the sequence $\{X_i,1\leq i\leq n\},1\leq k\leq n$, and $S_n=\sum_{i=1}^nX_i$.

Further, \cite{BBS13} includes a systematic study of properties of the class $\J$ such as its closure under weak tail equivalence with applications to random walks and infinitely divisible distributions.
Xu et al. \cite{XSWC14} construct some new heavy-tailed distributions in the class $\J$ to explain the differences between the class $\J$ and some well-established distribution classes, and describe the structure of the heavy-tailed distributions in the class $\J$. By contrast, using a simple transformation
which maps heavy-tailed to light-tailed distributions, we will illuminate the structure of the light-tailed distributions in $\J$ in the present paper. The structure is related to the convolution equivalent distribution class which, along with several heavy-tailed classes, will be
introduced below.

In this paper, unless otherwise stated, $F$ is a distribution with unbounded support in $[0,\infty)$ with tail $\overline F :=1-F$, $F*G$ is the convolution of two distributions $F$ and $G$, and $F^{*i}:=F* \dots *F$ is the $i$-fold convolution of $F$ with itself for $i\ge2$. In addition, all limits refer to $x$ tending to infinity.

Recall that a distribution $F$ belongs to the {\em subexponential} distribution class introduced by Chistyakov \cite{C64}, denoted by $F\in\mathcal{S}$, if
$$\overline {F^{*2}}(x)\sim2\overline F(x),$$
where the notation $g(x)\sim h(x)$ means that ratio $g(x)(h(x))^{-1}\rightarrow 1$ for two positive functions $g$ and $h$ supported on $[0,\infty)$. The class $\mathcal{S}$ is properly included in the following larger distribution class.

We say that a distribution $F$ belongs to the the long-tailed distribution class $\mathcal{L}$, if for all (or equivalently, for some) constant $t\neq 0$
$$\overline {F}(x+t)\sim\overline F(x).$$
The class $\mathcal{L}$ is an important subclass of the heavy-tailed distribution class $\mathcal{K}$ satisfying $\int_0^\infty e^{\lambda y}F(dy)=\infty$ for all $\lambda>0$. We will
denote its complement, the class of light-tailed distributions, by $\mathcal{K}^c$.

The following heavy-tailed distribution subclass was introduced by Feller \cite{F69}. A distribution $F$ belongs to the dominatedly-varying distribution class, denoted by $F\in\mathcal{D}$, if for all (or equivalently, for some) $t\in(0,1)$
$$\overline {F}(tx)=O(\overline F(x)),$$
where the notation $g(x)=O(h(x))$ means that $ \limsup g(x)(h(x))^{-1}<\infty$ for two positive functions $g$ and $h$ supported on $[0,\infty)$. The classes $\mathcal{L}$ and $\mathcal{D}$ are not comparable, see, for example,  \cite{EKM97}.

It is  proved in \cite{BBS13} that  $\mathcal{L}\cap\mathcal{J}=\mathcal{S}$ and $\mathcal{S}\cup\mathcal{D}\subset\mathcal{J}$. Xu et al. \cite{XSWC14} note that the last inclusion relationship is proper. In addition, the class $\J$ contains some common light-tailed distributions.

We say that a distribution $F$ belongs to the exponential distribution class, denoted by $F\in\mathcal{L}(\gamma)$ for some constant $\gamma>0$, if for all (or, equivalently, for some) $t \neq 0$
$$\overline {F}(x+t)\sim e^{-\gamma t}\overline F(x).$$
Further, if $F\in\mathcal{L}(\gamma)$ for some constant $\gamma>0,\ \int_0^\infty e^{\gamma y}F(dy)<\infty$ and $$\overline {F^{*2}}(x)\sim2\int_0^\infty e^{\gamma y}F(dy)\overline F(x),$$
then we say that $F$ belongs to the convolution equivalent distribution class, denoted by $F\in\mathcal{S}(\gamma)$ which was introduced in \cite{CNW73a},
\cite{CNW73b}.

Similarly to the heavy-tailed case, Proposition 5 of \cite{BBS13} shows that for any $\gamma>0$, we have $\mathcal{L}(\gamma)\cap\mathcal{J}=\mathcal{S}(\gamma)$ and hence $\mathcal{S}(\gamma)\subset\mathcal{J}$. Further, the authors state the following conjecture:
$$\J=(\J\cap\mathcal{K})\cup(\cup_{\beta>0}\mathcal{S}(\beta)).$$

In this paper, we give a negative answer to this conjecture, i.e.~we prove that the class $(\mathcal{J}\cap\mathcal{K}^c)\setminus(\cup_{\beta>0}\mathcal{S}(\beta))$ is non-empty. This fact indicates that the class $\J \cap \mathcal{K}^c$ -- just like $\J \cap \mathcal{K}$ -- cannot easily be expressed
in terms of well-established distribution classes. In particular, our results show that some light-tailed distributions without exponential nature
can still obey the principle of a single big jump.

Our main idea is to find a suitable transformation of distributions, which changes a heavy-tailed distribution into a light-tailed distribution belonging to the class $(\mathcal{J}\cap\mathcal{K}^c)\setminus(\cup_{\beta>0}\mathcal{S}(\beta))$. Some references use different methods to study the relationship between heavy-tailed distributions and light-tailed distributions. For example, Teugels \cite{T75}, Veraverbeke \cite{V77}, Embrechts and Goldie \cite{EG82}, Wang et al. \cite{WW11} use the Esscher transform of a distribution and Tang \cite{T08}  the product of random variables. Here, we use the following method, see Kl\"{u}ppelberg \cite{K89}. For a distribution $F$ and any fixed constant $\gamma>0$, we define the distribution $G$ as follows:
\begin{eqnarray}\label{102}
\overline{G}(x)=\textbf{1}(x<0)+\overline{F}(x)e^{-\gamma x}\textbf{1}(x\ge0),\ x\in(-\infty,\infty).
\end{eqnarray}

First, we want to find conditions on the distribution $F$, under which $G$ belongs to the class $\J$.
Further, we provide conditions on $F$ such that $G\in\J\cap\backslash(\cup_{\beta>0}\mathcal{S}(\beta))$.
To this end, we introduce further distribution classes.

We say that a distribution $F$ belongs to the generalized strongly subexponential distribution class $\mathcal{OS}^*$ introduced in \cite{WXCY14},
if its mean $\mu>0$ is finite and
$$C^\otimes(F):=\limsup\int_0^x\overline {F}(x-y)\overline F(y)dy(\overline F(x))^{-1}<\infty.$$
We have $C^\otimes(F)\ge2\mu$ for $F\in\mathcal{K}$ (\cite{FK07}). The class ${\cal{S}^*}:=\{F\in\mathcal{K},\   C^\otimes(F)=2\mu\}$
of {\em strongly subexponential} distributions was introduced in \cite{K89}. Some distributions with finite mean belonging to the class $\cal{L}\cap\cal{OS}^*\setminus\mathcal{S}^*$ can be found in Example 2.1 and Example 2.2 of \cite{WXCY14}. Clearly, the class $\mathcal{OS}^*$ can not contain the class $\cal{S}$. Conversely, the class $\cal{S}$ can not contain the class $\mathcal{OS}^*$, see the above-mentioned examples. The relation $\mathcal{S}^*\subset\mathcal{S}$ is proper, see \cite{KV91}, \cite{DFK04} and \cite{WXCY14}.

Similarly to the proper relationship $\mathcal{S}^*\subset\mathcal{S}$, the class $\mathcal{OS}^*$ is properly included in the generalized subexponential distribution class $\mathcal{OS}$ introduced in \cite{K90}, in which the distribution $F$ satisfies the following condition:
$$C^*(F):=\limsup\overline {F^{*2}}(x)(\overline {F}(x))^{-1}<\infty.$$
Proposition 5 (a) of \cite{BBS13} shows that $\mathcal{J}\subset\mathcal{OS}$. Some distributions in the class $\mathcal{L}\cap\mathcal{OS}\backslash S$ and hence in
$\mathcal{OS}\setminus\mathcal{J}$ can be found in \cite{L89}, \cite{LW12} and  \cite{WXCY14}.
On the other hand, the class $\mathcal{OS}$ is properly included in the following distribution class, see \cite{SW05}.

We say that a distribution $F$ belongs to the generalized long-tailed distribution class $\mathcal{OL}$, if for any $t>0$
$$C(F,t):=\limsup\overline {F}(x-t)(\overline {F}(x))^{-1}<\infty.$$
In addition, the relationship $\cup_{\beta\ge0}\mathcal{L}(\beta)\subset\mathcal{OL}$ is also proper, see for example, Example \ref{exam302} below, where the distribution $G_m\in\J\subset\mathcal{OS}\subset\mathcal{OL}$, but $G_m\notin\cup_{\beta>0}\mathcal{L}(\beta)$ for all $m\ge1$.
\vspace{0.2cm}

Now, we state the main results of the paper. First, we present two negative results by providing conditions on $F$ which imply $G\notin \J$ (with $G$ as defined
above). We will write $\J_f$ for the class of distributions in $\J$ with finite mean $\mu$ (and similarly for the classes $\mathcal{D}$ and $\S$). Note that we automatically have $\mu>0$ since distributions in $\J$ have unbounded support contained in $[0,\infty)$, as well as $\J\cap\mathcal{K}^c=\J_f\cap\mathcal{K}^c$.

\begin{pron}\label{pron101}
If $F$ belongs to the class $(\J_f\cap\mathcal{K})\backslash\mathcal{OS}^*$,
then $G\notin\mathcal{J}$. Further, the class $(\J_f\cap\mathcal{K})\backslash\mathcal{OS}^*$ is not empty.
\end{pron}

Even if the distribution $F$ belongs to the class $\J_f\cap\mathcal{K}\cap\mathcal{OS}^*$, however, the distribution $G$ does not necessarily belong to the class $\J$. In Example 2.1 of \cite{WXCY14}, there is a distribution $F\in(\mathcal{S}\cap\cal{OS}^*)\setminus\mathcal{S}^*$, thus $F\in\J_f$, but $G\notin\J$ according to following proposition.

\begin{pron}\label{pron102}
If $F$ belongs to the class $\S_f\setminus\mathcal{S}^*$, then $G\notin\mathcal{J}$.
\end{pron}

Next, we provide two positive results guaranteeing that $G\in\J$. According to Theorem 3.1 in  \cite{K89} and Proposition 5 in \cite{BBS13} we know that $G\in\mathcal{S}(\gamma)\subset\J$ under the condition that $F\in\mathcal{S}^*$. Therefore, if we require that $G\in(\J\cap\mathcal{K}^c)\setminus(\cup_{\beta>0}\mathcal{S}(\beta))$, we need to consider distributions $F$ outside the class $\mathcal{S}^*$.

\begin{pron}\label{pron103}
If $F$ belongs to the class ${\mathcal{D}}_f\setminus\mathcal{S}^*$, then $G\in(\mathcal{J}\cap\mathcal{K}^c)\setminus(\cup_{\beta>0}\mathcal{S}(\beta))$.
\end{pron}

Some examples of distributions in the class ${\mathcal{D}}_f\setminus\mathcal{S}$ can be found in \cite{EKM97}, \cite{WW11}, and in other references.
Obviously ${\mathcal{D}}_f\subset\mathcal{OS}^*$. The following proposition shows that there also exist distributions $F$ outside the class $\mathcal{D}$ such that $G\in(\mathcal{J}\cap\mathcal{K}^c)\setminus(\cup_{\beta>0}\mathcal{S}(\beta))$.

\begin{pron}\label{pron104}
There exists a distribution $F\in(\J_f\cap\mathcal{K})\setminus(\mathcal{S}\cup\mathcal{D})$ such that
$G\in \mathcal{J}$. Any such $F$ is automatically in $\mathcal{OS}^*$ and $G\in(\mathcal{J}\cap\mathcal{K}^c)\setminus(\cup_{\beta>0}\mathcal{S}(\beta))$.
\end{pron}

Based on Proposition \ref{pron104}, we provide additional information on the structure of the class $\J\cap\mathcal{K}^c$. Combined with Theorem 1.2 of
\cite{XSWC14}, we are thus able to explain the structure of the class $\J$. To this end, we define the distribution classes
$$\J(\mathcal{K}^c)_1:=\J\cap\mathcal{K}^c\cap\{F:\mbox{ there exists some } F_1\in\mathcal{S}(\gamma) \mbox{ for some } \gamma>0\mbox{ such that } \overline{F}(x)\approx\overline{F_1}(x) \}$$
and $\J(\mathcal{K}^c)_2:=(\J\cap\mathcal{K}^c)\backslash\J(\mathcal{K}^c)_1$, where the notation $g(x)\approx h(x)$ means $g(x)=O(h(x))$ and $h(x)=O(g(x))$ for two positive functions supported on $[0,\infty)$. We say that a distribution $F$ is {\em weakly tail equivalent} to another distribution $F_1$, if $\overline{F}(x)\approx\overline{F_1}(x)$.

Correspondingly, the following two heavy-tailed distribution classes are introduced by \cite{XSWC14}:
$$\J(\mathcal{K})_1:=\J\cap\mathcal{K}\cap\{F:\mbox{ there exists some } F_1\in\mathcal{S}\mbox{ such that } \overline{F}(x)\approx\overline{F_1}(x) \}$$
and $\J(\mathcal{K})_2:=(\J\cap\mathcal{K}^c)\backslash\J(\mathcal{K})_1$. Further, we define the distribution classes
$$\J_1:=\J(\mathcal{K})_1\cup\J(\mathcal{K}^c)_1 \mbox{ and } \J_2:=\J(\mathcal{K})_2\cup\J(\mathcal{K}^c)_2.$$
In addition, we write $\mathcal{S}(0):=\mathcal{S}$. Then we have the following result.

\begin{thm}\label{thm101}
i) $\J\cap\mathcal{K}^c=\J(\mathcal{K}^c)_1\cup\J(\mathcal{K}^c)_2$, and $\J(\mathcal{K}^c)_i\backslash(\cup_{\beta>0}\mathcal{S}(\beta))\neq\emptyset,\ i=1,2$.

ii)  $\J=\J_1\cup\J_2$, and $\J_i\backslash(\cup_{\beta\ge0}\mathcal{S}(\beta))\neq\emptyset,\ i=1,2$.
\end{thm}

\begin{remark}\label{remark101}
In the above results, the condition  $F\in\J_f$, that is $0<\mu<\infty$, is indispensable in certain sense, because if the mean of $F$ is infinite,
then formula \eqref{206} together with \eqref{new} show that $G\notin\mathcal{OS}$, thus $G\notin\mathcal{J}$.
\end{remark}

We will prove Propositions \ref{pron101}-\ref{pron104} and Theorem \ref{thm101} in Section 3. To this end, we give three lemmas in the next section.

\section{Some lemmas}
\setcounter{thm}{0}\setcounter{Corol}{0}\setcounter{lemma}{0}\setcounter{pron}{0}
\setcounter{remark}{0}\setcounter{exam}{0}\setcounter{equation}{0}

First, we state the following simple lemma.

\begin{lemma}\label{lemma202}
Let $F$ be a distribution with finite and positive mean $\mu$. Then the following three statements are mutually equivalent: $F\in\mathcal{OS}^*$,  $G\in\mathcal{OS}$ and $G\in\mathcal{OS}^*$.
\end{lemma}

\proof According to (\ref{202}), we have
\begin{eqnarray}\label{206}
\frac{\overline {G^{*2}}(x)}{\overline {G}(x)}=\frac{\overline {F^{*2}}(x)}{\overline {F}(x)}+\frac{\gamma\int_{0}^{x}\overline {F}(x-y)\overline {F}(y)dy}{\overline {F}(x)}.
\end{eqnarray}
Since $\mathcal{OS}^* \subset \mathcal{OS}$, \eqref{206} shows that $G\in\mathcal{OS}$ iff $F\in\mathcal{OS}^*$. Further, from the definitions of $\mathcal{OS}^*$ and of $G$ we see that
$F\in\mathcal{OS}^*$ iff $G\in\mathcal{OS}^*$, so the statement of the lemma follows.\hfill$\Box$\\

Next, we provide a necessary and sufficient condition on $F \in \J_f$ for its associated distribution $G$ to be in $\J$ (with $G$ defined as in \eqref{102}).

\begin{lemma}\label{lemma201}
If the distribution $F$ belongs to the class $\J_f$, then $G\in\mathcal{J}$ if and only if
\begin{eqnarray}\label{201}
\lim_{K\to\infty}\liminf T^{(X)}(x;K):=\lim_{K\to\infty}\liminf
\frac{2\int_{0}^{K}\overline {F}(x-y)\overline {F}(y)dy}
{\int_{0}^{x}\overline {F}(x-y)\overline {F}(y)dy}=1.
\end{eqnarray}
\end{lemma}

\proof~$\Leftarrow$ Let $Y_1$ and $Y_2$ be two independent random variables with common distribution $G$. Denote $S_2^{(Y)}:=Y_1+Y_2$, $Y_{2,2}:=\min\{Y_1,Y_2\}$ and for any constant $K>0$ and $x>2K$
\begin{eqnarray*}
B^{(Y)}(x;K):=\P(Y_{2,2}\leq K|S_2^{(Y)}> x)=2(\overline {G^{*2}}(x))^{-1}\int_{0}^{K}\overline {G}(x-y)dG(y).
\end{eqnarray*}
an easy calculation shows that
\begin{eqnarray}\label{202}
\overline {G^{*2}}(x)=\Big(\overline {F^{*2}}(x)+2\gamma\int_{\frac{x}{2}}^{x}\overline {F}(x-y)\overline {F}(y)dy\Big)e^{-\gamma x}
\end{eqnarray}
and
$$2\int_{0}^{K}\overline {G}(x-y)dG(y)=\Big(2\int_{0}^{K}\overline {F}(x-y)F(dy)+2\gamma\int_{0}^{K}\overline {F}(x-y)\overline {F}(y)dy\Big)e^{-\gamma x}.$$
According to the  elementary inequality $\frac{a+c}{b+c}\geq\frac{a}{b}$ for any $c>0$ and any $b\geq a>0$, and noting  that
$$\int_{\frac{x}{2}}^{x}\overline {F}(x-y)\overline {F}(y)dy=\int_{0}^{\frac{x}{2}}\overline {F}(x-y)\overline {F}(y)dy
\ge\int_{0}^{K}\overline {F}(x-y)\overline {F}(y)dy$$
for $x\ge 2K$, we have
\begin{eqnarray*}
&&B^{(Y)}(x;K)=\Big(\overline {F^{*2}}(x)+2\gamma\int_{\frac{x}{2}}^{x}\overline {F}(x-y)\overline {F}(y)dy\Big)^{-1}\nonumber\\
& &\ \ \ \cdot\Bigg(\Big(2\int_{0}^{K}\overline {F}(x-y)F(dy)+2\gamma\int_{\frac{x}{2}}^{x}\overline {F}(x-y)\overline {F}(y)dy\Big)-\Big(\overline {F^{*2}}(x)+2\gamma\int_{\frac{x}{2}}^{x}\overline {F}(x-y)\overline {F}(y)dy\Big)\nonumber\\
&&\ \ \ \ \ \
\ \ \ +\Big(\overline {F^{*2}}(x)+2\gamma\int_{0}^{K}\overline {F}(x-y)\overline {F}(y)dy\Big)\Bigg)\nonumber\\
&\geq&\lo(\overline {F^{*2}}(x)\ro)^{-1}\Big(2\int_{0}^{K}\overline {F}(x-y)F(dy)\Big)+\Big(\int_{\frac{x}{2}}^{x}\overline {F}(x-y)\overline {F}(y)dy\Big)^{-1}
\Big(\int_{0}^{K}\overline {F}(x-y)\overline {F}(y)dy\Big)-1.
\end{eqnarray*}
Hence, by $F\in\mathcal{J}$ and (\ref{201}), we obtain $\lim_{K\rightarrow\infty}\liminf B^{(Y)}(x;K)\to 1$, that is $G\in\mathcal{J}$.

$\Rightarrow$ Since $F\in \cal{J}$, we have $F\in \cal{OS}$. Because
\begin{equation}\label{new}
\int_{\frac{x}{2}}^{x}\overline {F}(x-y)\overline {F}(y)dy\ge\overline {F}(x)\int_0^{\frac{x}{2}}\overline {F}(y)dy\sim\overline {F}(x)\mu,
\end{equation}
we know that there exists a constant $C>0$ such that $\overline {F^{*2}}(x)\leq C\int_{\frac{x}{2}}^{x}\overline {F}(x-y)\overline {F}(y)dy$ for $x$ large enough. Therefore,
for $x$ large enough,

\begin{eqnarray}\label{205}
&&B^{(Y)}(x;K)=\Big(\overline {F^{*2}}(x)+2\gamma\int_{\frac{x}{2}}^{x}\overline {F}(x-y)\overline {F}(y)dy\Big)^{-1}\nonumber\\
& &\ \ \  \  \ \ \ \ \cdot\bigg(\Big(2\int_{0}^{K}\overline {F}(x-y)F(dy)+2\gamma\int_{\frac{x}{2}}^{x}\overline {F}(x-y)\overline {F}(y)dy\Big)\nonumber\\
& &\ \ \ \ \ \ \ \ \ \ \ \ \ \ \ \ \ -\Big(2\gamma\int_{\frac{x}{2}}^{x}\overline {F}(x-y)\overline {F}(y)dy-2\gamma\int_{0}^{K}\overline {F}(x-y)\overline {F}(y)dy\Big)\bigg)\nonumber\\
&\leq&1-\Big((C+2\gamma)\int_{\frac{x}{2}}^{x}\overline {F}(x-y)\overline {F}(y)dy\Big)^{-1}\Big(2\gamma\int_{\frac{x}{2}}^{x}\overline {F}(x-y)\overline {F}(y)dy-2\gamma\int_{0}^{K}\overline {F}(x-y)\overline {F}(y)dy\Big)
\nonumber\\
&=&\frac{C}{C+2\gamma}+\frac{2\gamma}{C+2\gamma}\Big(\int_{\frac{x}{2}}^{x}\overline {F}(x-y)\overline {F}(y)dy\Big)^{-1}\int_{0}^{K}\overline {F}(x-y)\overline {F}(y)dy\nonumber\\
&=&\frac{C}{C+2\gamma}+\frac{2\gamma}{C+2\gamma}T^{(X)}(x;K).
\end{eqnarray}

By $G\in \cal{J}$ and (\ref{205}) we have
\begin{eqnarray*}
1\ge\lim_{K\to\infty}\liminf_{x\to\infty}T^{(X)}(x;K)
\ge&\lim_{K\to\infty}\liminf_{x\to\infty}\Big(B^{(Y)}(x;K)-\frac{C}{C+2\gamma}\Big)
\frac{C+2\gamma}{2\gamma}=1,
\end{eqnarray*}
that is (\ref{201}) holds.\hfill$\Box$
\vspace{0.2cm}

Finally, we state a modified version of Lemma 4.1 in \cite{XSWC14}. For this, we first introduce some related function classes.

Similarly to the distribution class $\mathcal{L}$ and the distribution class $\mathcal{OL}$, we say that a positive function $f$ supported on $[0,\infty)$ belongs to the long-tailed function class $\mathcal{L}_d$, if for all $t\neq 0$
$$f(x-t)\sim f(x);$$
and say that the function $f$ belongs to the generalized long-tailed function class $\mathcal{OL}_d$, if for all $t\neq 0$
$$c(f,t):=\limsup f(x-t)(f(x))^{-1}<\infty.$$
Further, we say that a function $f$ belongs to the
strongly subexponential function class $\mathcal{S}_d$, if $f\in\mathcal{L}_d$, $\int_0^\infty f(y)dy<\infty$ and
$$\int_0^x f(x-y)f(y)dy\sim2f(x)\int_0^\infty f(y)dy.$$
Obviously, a distribution $F\in\cal{L}$ if and only if $f:=\overline{F}\in\mathcal{L}_d$, $F\in\cal{OL}$ if and only if $\overline{F}\in\mathcal{OL}_d$ and $F\in\mathcal{S}^*$ if and only if $\overline{F}\in\mathcal{S}_d$.

\begin{lemma}\label{lemma203}
Assume that the  function $f\in\mathcal{OL}_d$ satisfies
\begin{eqnarray}
\limsup_{t\rightarrow\infty} c(f,t)=\infty.\label{207}
\end{eqnarray}
Then $f$ is not weakly equivalent to any long-tailed function.
\end{lemma}

\noindent{\bf Proof}\ \ We assume there exists a function $f_1\in\mathcal{L}_d$ and $f(x)\approx f_1(x)$. Then there are two constants $0<c_1\leq c_2<\infty$ such that
\begin{eqnarray}\label{208}
c_1= \liminf(f_1(x))^{-1}f(x)\leq \limsup(f_1(x))^{-1}f(x)= c_2.
\end{eqnarray}
By $f_1\in\mathcal{L}_d$ and (\ref{208}) for any $0<t<\infty$ we have
\begin{eqnarray}\label{209}
(f(x))^{-1}f(x-t)
\lesssim(c_1f_1(x))^{-1}c_2f_1(x-t)\sim&c_1^{-1}c_2,
\end{eqnarray}
where $f(x)\lesssim g(x)$ means $\limsup g(x)(h(x))^{-1}\leq1$ for two positive functions $g$ and $h$ supported on $[0,\infty)$. Obviously, (\ref{209}) contradicts  (\ref{207}), hence the conclusion of the lemma holds.\hfill$\Box$

\section{Proofs}
\setcounter{thm}{0}\setcounter{Corol}{0}\setcounter{lemma}{0}\setcounter{pron}{0}
\setcounter{remark}{0}\setcounter{exam}{0}\setcounter{equation}{0}

\subsection{Proof of Proposition \ref{pron101}}

If $F$ belongs to the class $(\J_f\cap\mathcal{K})\backslash\mathcal{OS}^*$, then,
according to Lemma \ref{lemma202}, $G \notin \mathcal{OS}$ and hence $G\notin{\J}$.

Now we prove that $(\J_f\cap\mathcal{K})\backslash\mathcal{OS}^*\neq\emptyset$ by the following example from Section 3.8 in \cite{FKZ13}
which the authors constructed to prove that the integrated tail distribution of $F \in \S_f$ is not necessarily in the class $\S$.
\begin{exam}\label{exam300}
Define $a_0=0$, $a_1=1$ and $a_{n+1}=\frac{e^{a_n}}{a_n}$. Since $\frac{e^x}{x}$ is increasing for $x\geq1$,
the sequence $a_n$ is increasing to infinity and $a_n=o(a_{n+1})$, as $n\to\infty$. Define the  distribution $F$ as follows:
\begin{eqnarray}\label{3002}
\overline{F}(x)=\textbf{\emph{1}}(x<0)+\sum\limits_{n=0}^{\infty}e^{-\lo(a_n+\frac{(x-a_n^2)}{a_n+a_{n+1}}\ro)}
\textbf{\emph{1}}(a_n^2\le x<a_{n+1}^2).
\end{eqnarray}
Foss et al. \cite{FKZ13} proved that $F\in\S_f\setminus{\cal{S^{*}}}$. We show that $F\notin\cal{OS^*}$.
\end{exam}
By (\ref{3002}), we know that
\begin{eqnarray*}
(\overline{F}(a_{n+1}^2))^{-1}\int_{0}^{a_{n+1}^2}\overline{F}(y)\overline{F}(a_{n+1}^2-y)dy
&\geq& e^{a_{n+1}}\int_{a_{n+1}^2/2}^{a_{n+1}^2}\overline{F}(y)\overline{F}(a_{n+1}^2-y)dy\nonumber\\
&=&\int_{0}^{a_{n+1}^2/2}\overline{F}(y)e^{\frac{y}{a_n+a_{n+1}}}dy\nonumber\\
&\geq&\int_{a_{n}^2}^{a_{n+1}^2/2}\overline{F}(y)e^{\frac{y}{a_n+a_{n+1}}}dy
\nonumber\\
&=&\int_{a_{n}^2}^{a_{n+1}^2/2}e^{-\lo(a_n-\frac{a_n^2}{a_n+a_{n+1}}\ro)}dy\nonumber\\
&\geq&(2^{-1}a_{n+1}^2-a_{n}^2)e^{-a_n}\to\infty,~~~n\rightarrow\infty.
\end{eqnarray*}
Thus $F\notin\cal{OS^*}$.


\subsection{Proof of Proposition \ref{pron102}}

If $F\in \S_f\setminus\cal{S}^*$, then by Theorem 2.1 in \cite{K89}, $G\in \cal{L}(\gamma)\setminus\cal{S}(\gamma)$ thus, by Proposition 5 in \cite{BBS13},
$G\notin \cal{J}$.

\subsection{Proof of Proposition \ref{pron103}}

Obviously, $G\in\mathcal{K}^c$. In the following, we prove that $G\in\mathcal{J}$. To this end, we only need to prove that (\ref{201}) holds according to Lemma \ref{lemma201} and the fact that $\mathcal{D}\subset \mathcal{J}$. Since $\int_0^\infty\overline F(y)dy<\infty$ and $F\in\mathcal{D}$, that is $D(F):=\sup_{x\ge0}\frac{\overline {F}(\frac{x}{2})}{\overline {F}(x)}<\infty$, we obtain that,
for $x \ge 2K$,

\begin{eqnarray*}
1&\ge&\frac{\int_{0}^{K}\overline {F}(x-y)\overline {F}(y)dy}{\int_{\frac{x}{2}}^{x}\overline {F}(x-y)\overline {F}(y)dy}\nonumber\\
&=&1-\frac{\int_{K}^{\frac{x}{2}}\overline {F}(x-y)\overline {F}(y)dy}{\int_{0}^{\frac{x}{2}}\overline {F}(x-y)\overline {F}(y)dy}\nonumber\\
&\geq&1-\frac{\overline {F}(\frac{x}{2})\int_{K}^{\frac{x}{2}}\overline {F}(y)dy}{\overline {F}(x)\int_{0}^{\frac{x}{2}}\overline {F}(y)dy}\nonumber\\
&\geq&1-\frac{D(F)\int_{K}^{\frac{x}{2}}\overline {F}(y)dy}{\int_{0}^{\frac{x}{2}}\overline {F}(y)dy}\nonumber\\
&\to&1,~~~x\rightarrow\infty,\quad K\rightarrow\infty.
\end{eqnarray*}
Hence (\ref{201}) holds.

Finally, from Theorem 2.1 in \cite{K89}, we know that $G\in \cal{S}(\gamma)$ if and only if $F\in \cal{S^*}$, and since
$F\notin \cal{S^*}$ we obtain  $G\notin \cal{S}(\gamma)$. Now take any $\beta>0$ which is different from $\gamma$ and assume that
$G\in \cal{S}(\beta)$. Then $G\in \cal{L}(\beta)$ which, together with the definition of $G$, implies
$$
e^{-\gamma t} \overline {F}  (x+t) \sim e^{-\beta t} \overline {F}  (x),
$$
which is impossible in case $\beta < \gamma$ and implies $F \in \cal{L}(\beta - \gamma)$ in case  $\beta > \gamma$ contradicting
the assumption that $F$ is in $\cal{D}$ and therefore heavy-tailed.


\subsection{Proof of Proposition \ref{pron104}}

Let us first show the last statements in Proposition  \ref{pron104}. The fact that $F \in \mathcal{OS}^*$ follows from Lemma \ref{lemma202} since $G \in \J \subset \mathcal{OS}$
and $F \notin \S$ implies $F \notin \S^*$ from which we get $G \notin \cup_{\beta>0}\mathcal{S}(\beta)$ as in the last part of the proof of Proposition \ref{pron103}.

Proposition \ref{pron104} therefore follows from the following example which is Example 2.1 in \cite{XSWC14}.

\begin{exam}\label{exam301}
Assume that $F_1\in\mathcal{S}^*\cap\mathcal{K}_1$ is a continuous distribution and let $y_0\ge0$ and $a>1$ be two constants such that $a\overline{F}_1(y_0)\le1$, where $\mathcal{K}_1=\{H\in\mathcal{K},\ \int_0^\infty y^b H(dy)<\infty\ for\ all\ b>0\}$. For example,
we can take $\overline{F}_1(x)=\textbf{\emph{1}}(x<0)+e^{-\sqrt{x}}\textbf{\emph{1}}(x\ge0),\ x\in(-\infty,\infty)$ and $y_0=(\ln a)^2$. Let $F$ be a distribution such that for all $x\in(-\infty,\infty)$,

\begin{eqnarray}\label{exam-e101first}
\overline{F}(x)&=&\overline{F}_1(x)\textbf{\emph{1}}(x<x_1)+
\sum_{i=1}^{\infty}\Big(\overline{F}_1(x_i)\textbf{\emph{1}}(x_i\le x<y_i)+\overline{F}_1(x)\textbf{\emph{1}}(y_i\le x<x_{i+1})\Big),
\end{eqnarray}
where $\{x_i,i\ge1\}$ and $\{y_i, i\ge 1\}$ are two sequences of positive constants satisfying $x_i<y_i<x_{i+1}$ and $\overline{F}_1(x_i)=a\overline{F}_1(y_i),\ i\ge1$. Obviously, $\overline{F}(x)\approx\overline{F_1}(x)$. Xu et al. \cite{XSWC14} proved that $F\in\J_f\setminus(\cal{L}\cup\cal{D})$.
\end{exam}
Let $\overline{G}_1(x)=\overline{F}_1(x)e^{-\gamma x},~x\in(-\infty,\infty)$, then $G_1\in\mathcal{S}(\gamma)\subset\J$ and $\overline{G}(x)\approx\overline{G_1}(x)$, thus $G\in\J$.
\hfill$\Box$

\subsection{Proof of Theorem \ref{thm101}}

Note that part ii) of the theorem is a simple consequence of part i) together with Theorem 1.2 in \cite{XSWC14}.
Note further that Example \ref{exam301} in the previous subsection proves one half of Theorem \ref{thm101} i), namely that $\J(\mathcal{K}^c)_1\backslash \cup_{\beta >0}\mathcal{S}(\beta)\neq \emptyset$. The other half
follows from the next example which is another example satisfying the claim in Proposition  \ref{pron104} but in this case $G$ is not
weakly tail equivalent to any convolution equivalent distribution.

\begin{exam}\label{exam302}
For any a positive integer $m$, choose any constant $\alpha\in
(2+3m^{-1},\infty)$ and any constant $x_1>4^{\alpha}$. For all
integers $n\geq1$, let $x_{n+1}=x_n^{1+{\alpha}^{-1}}$. Clearly,
$x_{n+1}>4x_n$ and $x_n\to\infty$ as $n\to\infty$. Now, define the
distribution $F$ as follows£º
\begin{eqnarray}\label{exam-e101}
&&\overline{F}(x)=\emph{\textbf{1}}(x<0)+(x_1^{-1}(x_1^{-\alpha}-1)x+1)\emph{\textbf{1}}(0\leq x< x_1)\nonumber\\
&+&\sum\limits_{n=1}^{\infty}((x_n^{-\alpha}+(x_n^{-\alpha-2}-x_n^{-\alpha-1})(x-x_n))\emph{\textbf{1}(}x_n\leq x<2x_n)
+x_n^{-\alpha-1}\emph{\textbf{1}}(2x_n\leq x< x_{n+1})).
\end{eqnarray}
Define the distribution $F_m$ by $\overline{F}_m(x)=(\overline{F}(x))^m$ for all $x\in(-\infty,\infty)$. It was shown in  \cite{XSWC14} that
$F_m\in(\J_f\cap\mathcal{K})\setminus(\cal{L}\cup\mathcal{D})$. Further, for any fixed $\gamma>0$ let
$$\overline{G}_m(x)=\overline{F}_m(x)e^{-\gamma x},~x\in(-\infty,\infty),$$
then we prove that $G_m\in(\mathcal{J}\cap\mathcal{K}^c)\setminus(\cup_{\beta>0}\cal{S}(\beta))$ and $G_m$ is not weakly tail equivalent to any convolution equivalent distribution.
\end{exam}

According to Proposition 12 in \cite{BBS13}, we only need to prove the conclusion in case $m=1$.

 From (\ref{exam-e101}), it is easy to see that when
$x\geq x_1,$
\begin{eqnarray}
x^{-\alpha-1}\leq\overline{F}(x)\leq 2^{\alpha}x^{-\alpha}.\label{e101-2}
\end{eqnarray}
Hence, $F$ is heavy-tailed and $G:=G_1$ is light-tailed. Obviously, $F$ has a finite mean which we denote by $\mu$. In fact, we even have
\begin{eqnarray}
\int_{0}^{\infty}y^4\overline{F}(y)dy<\infty.\label{e101-3}
\end{eqnarray}

First, by $F\notin{\cal{L}}$, we have $G\notin{\cal{L}(\gamma)}$ and -- as in the proof of Proposition \ref{pron103} --
$G \notin \cup_{\beta>0}\cal{S}(\beta))$.

Next we prove $G\in{\cal{J}}.$ By Lemma \ref{lemma201} we only need to prove that
(\ref{201}) holds. For this end, we estimate $T^{(X)}(x;K)$ in each of the five
cases $x_n\leq x< x_n+K$, $x_n+K\leq x< \frac{3}{2}x_n$, $\frac{3}{2}x_n\leq x< 2x_n$,
$2x_n\leq x< 2x_n+K$ and $2x_n+K\leq x< x_{n+1}$, respectively.

When $x\in[x_n,x_n+K)$ and $0\leq y \leq K$, $2x_{n-1}\le \frac{x_n}{2}\le \frac{x}{2}\le x_n-K\leq x-y \leq x_n$ for large enough $n$. Thus, by (\ref{exam-e101}), we have $\overline {F}(\frac{x}{2})=\overline {F}(x_n)$ and
\begin{eqnarray}\label{e101-5}
T^{(X)}(x;K)&\ge&\Big(\overline {F}(\frac{x}{2})\int_{0}^{\frac{x}{2}}\overline {F}(y)dy\Big)^{-1}
\overline {F}(x_n)\int_{0}^{K}\overline {F}(y)dy\nonumber\\
&\geq&\Big(\int_{0}^{\infty}\overline {F}(y)dy\Big)^{-1}\int_{0}^{K}\overline {F}(y)dy\nonumber\\
&\to&1,~~~K\rightarrow\infty.
\end{eqnarray}

When $x\in[x_n+K,\frac{3}{2}x_n)$, $2x_{n-1}\le \frac{x_n+K}{2}\le \frac{x}{2}\leq \frac{3}{4}x_n\le x_n$ for large enough $n$. Thus, by (\ref{exam-e101}) and (\ref{e101-3}), we have

\begin{eqnarray}
&&T^{(X)}(x;K)\geq\Big(\overline {F}(\frac{x}{2})\int_{\frac{x}{2}}^{x_n}\overline {F}(x-y)dy+\int_{x_n}^{x}\overline {F}(x-y)\overline {F}(y)dy\Big)^{-1}\overline {F}(x)\int_{0}^{K}\overline {F}(y)dy\nonumber\\
&\geq&\Bigg(x_n^{-\alpha}\int_{x-x_n}^{\frac{x}{2}}\overline {F}(y)dy+\int_{0}^{x-x_n}\Big(\overline {F}(x)+x_n^{-\alpha-1}y\Big)\overline {F}(y)dy\Bigg)^{-1}\overline {F}(x)\int_{0}^{K}\overline {F}(y)dy\nonumber\\
&\geq&\lo(\frac{x_n^{-\alpha}}{\overline {F}(\frac{3}{2}x_n)}\int_{K}^{\frac{3}{4}x_n}\overline {F}(y)dy+\int_{0}^{\frac{x_n}{2}}\Big(1+\frac{x_n^{-\alpha-1}y}{\overline {F}(\frac{3}{2}x_n)}\Big)\overline {F}(y)dy\ro)^{-1}\int_{0}^{K}\overline {F}(y)dy\nonumber\\
&\gtrsim&\Big(2\int_{K}^{\frac{3}{4}x_n}\overline {F}(y)dy+\int_{0}^{\frac{x_n}{2}}\overline {F}(y)dy\Big)^{-1}\int_{0}^{K}\overline {F}(y)dy\nonumber\\
&\ge&\Big(2\int_{K}^{\infty}\overline {F}(y)dy+\mu\Big)^{-1}\int_{0}^{K}\overline {F}(y)dy\nonumber\\
&\to&1,~~~K\rightarrow\infty.
\label{e101-6}
\end{eqnarray}

When $x\in[\frac{3}{2}x_n,2x_n)$ and $0\leq y \leq K$, $\frac{3}{2}x_n-K\leq x-y \leq 2x_n-K$, $n\geq2$. Thus, by (\ref{exam-e101}) and (\ref{e101-3}), we have

\begin{eqnarray}
&&T^{(X)}(x;K)=\lo(\Big(\int_{\frac{x}{2}}^{x_n}+\int_{x_n}^{x}\Big)\overline {F}(x-y)\overline {F}(y)dy\ro)^{-1}\int_{0}^{K}\Big(\overline {F}(x)+(x_n^{-\alpha-1}-x_n^{-\alpha-2})y\Big)\overline {F}(y)dy\nonumber\\
&\gtrsim&\Big(x_n^{-\alpha}\int_{\frac{x_n}{2}}^{x_n}\overline {F}(y)dy+\int_{0}^{x-x_n}(\overline {F}(x)+x_n^{-\alpha-1}y)\overline {F}(y)dy\Big)^{-1}\Big(\int_{0}^{K}(\overline {F}(x)+x_n^{-\alpha-1}y)\overline {F}(y)dy\Big)\nonumber\\
&\geq&1-\Big(\int_{0}^{x-x_n}(\overline {F}(x)+x_n^{-\alpha-1}y)\overline {F}(y)dy\Big)^{-1}\Big(\int_{K}^{x-x_n}(\overline {F}(x)+x_n^{-\alpha-1}y)\overline {F}(y)dy+\frac{1}{2}x_n^{1-2\alpha}\Big)\nonumber\\
&\geq&1-\Big(\int_{0}^{\frac{x_n}{2}}\overline {F}(y)dy\Big)^{-1}\Big(\int_{K}^{x_n}\overline {F}(y)dy+\frac{1}{2}x_n^{2-\alpha}\Big)-\Big(\int_{0}^{\frac{x_n}{2}}y\overline {F}(y)dy\Big)^{-1}\int_{K}^{x_n}y\overline {F}(y)dy\nonumber\\
&\sim&1-\mu^{-1}\int_{K}^{\infty}\overline {F}(y)dy-\Big(\int_{0}^{\infty}y\overline {F}(y)dy\Big)^{-1}\Big(\int_{K}^{\infty}y\overline {F}(y)dy\Big)\nonumber\\
&\to&1,~~~K\rightarrow\infty.
\label{e101-7}
\end{eqnarray}

When $x\in[2x_n,2x_n+K)$ and $0\leq y \leq K$, $2x_n\le 2x_n+K-y\le x_{n+1}$ such that $\overline {F}(2x_n+K-y)=\overline {F}(2x_n)=x_n^{-\alpha-1}$; when $x\in[2x_n,2x_n+K)$ and $K\leq y \leq x_n$, $x_n\le2x_n-y\le2x_n$; and when $x\in[2x_n,2x_n+K)$ and $x_n\leq y \leq x_n+\frac{K}{2}$, $x_n-\frac{K}{2}\le2x_n-y\le x_n$ such that $\overline {F}(2x_n-y)=\overline {F}(x_n)=x_n^{-\alpha}$, $n\geq2$. Thus, by (\ref{exam-e101}) and (\ref{e101-3}), we have

\begin{eqnarray}
&&T^{(X)}(x;K)=1-\Big(\int_{0}^{\frac{x}{2}}\overline {F}(x-y)\overline {F}(y)dy\Big)^{-1}\int_{K}^{\frac{x}{2}}\overline {F}(x-y)\overline {F}(y)dy\nonumber\\
&\geq&1-\Big(\int_{0}^{K}\overline {F}(2x_n+K-y)\overline {F}(y)dy\Big)^{-1}
\Big(\int_{K}^{x_n}+\int_{x_n}^{x_n+\frac{K}{2}}\Big)\overline {F}(2x_n-y)\overline {F}(y)dy\nonumber\\
&\geq&1-\Big(\int_{0}^{K}x_n^{-\alpha-1}\overline {F}(y)dy\Big)^{-1}\Big(\int_{K}^{x_n}x_n^{-\alpha-1}(1+y)\overline {F}(y)dy+\frac{K}{2}x_n^{-2\alpha}\Big)\nonumber\\
&\gtrsim&1-\Big(\int_{0}^{K}\overline {F}(y)dy\Big)^{-1}\int_{K}^{\infty}(1+y)\overline {F}(y)dy\nonumber\\
&\to&1,~~~K\rightarrow\infty.
\label{e101-8}
\end{eqnarray}

When $x\in[2x_n+K,x_{n+1})$ and $0\leq y \leq K$, $2x_n\leq x-y \leq x_{n+1}-K$ such that $\overline {F}(x-y)=\overline {F}(2x_n)=x_n^{-\alpha-1}$; when $x\in[2x_n+K,x_{n+1})$ and $0\leq y \leq x-2x_n$, $2x_n\le x-y\le x_{n+1}$ such that $\overline {F}(x-y)=x_n^{-\alpha-1}$; and when $K\leq y \le x_n+\frac{K}{2}$, $x_n\le2x_n+K-y\le 2x_n$ $n\geq2$. Thus, by (\ref{exam-e101}) and (\ref{e101-3}), we have

\begin{eqnarray}
&&T^{(X)}(x;K)=\lo(\Big(\int_{\frac{x}{2}}^{2x_n}+\int_{2x_n}^{x}\Big)\overline {F}(x-y)\overline {F}(y)dy\ro)^{-1}x_n^{-\alpha-1}\int_{0}^{K}\overline {F}(y)dy\nonumber\\
&\geq&\Big(\int_{x_n+\frac{K}{2}}^{2x_n}\overline {F}(2x_n+K-y)\overline {F}(y)dy+
x_n^{-\alpha-1}\int_{0}^{x-2x_n}\overline {F}(y)dy\Big)^{-1}x_n^{-\alpha-1}\int_{0}^{K}\overline {F}(y)dy\nonumber\\
&\geq&\Big(\int_{K}^{x_n+\frac{K}{2}}(1+y)\overline {F}(y)dy+
\int_{0}^{x_{n+1}-2x_n}\overline {F}(y)dy\Big)^{-1}\int_{0}^{K}\overline {F}(y)dy\nonumber\\
&\sim&\Big(\int_{K}^{\infty}(1+y)\overline {F}(y)dy+
\mu\Big)^{-1}\int_{0}^{K}\overline {F}(y)dy\nonumber\\
&\to&1,~~~K\rightarrow\infty.
\label{e101-9}
\end{eqnarray}

According to (\ref{e101-5})-(\ref{e101-9}), we get $G\in \cal{J}$.

It remains to show that $G$ is not weakly tail equivalent to any  convolution equivalent distribution. It is clear that $G$ can not be weakly tail equivalent to any $H \in \mathcal{S}(\beta)$ in case $\beta \neq \gamma$.
Assume that there exists a distribution $H\in \mathcal{S}(\gamma)$ such that $\overline{G}(x)\approx \overline{H}(x)$. Write $f_1(x)=e^{\gamma x}\overline{H}(x),\ x\in(-\infty,\infty)$, then $f_1\in\mathcal{S}_d\subset\mathcal{L}_d$ by Theorem 2.1 in \cite{K89} and $f_1(x)\approx \overline{F}(x)$. However, according Lemma \ref{lemma203} and
$$(\overline{F}(2x_n))^{-1}\overline{F}(2x_n-t)=1+t-tx_n^{-1}\sim1+t,~~~n\to\infty,$$
$F$ is not weakly tail equivalent to any long-tailed function, which contradicts $\overline{F}(x)\approx f_1(x)$ and therefore the proof of Theorem \ref{thm101} is complete.\hfill$\Box$


\begin{thebibliography}{99}
\bibitem{DFK04} D. Denisov, S. Foss and  D. Korshunov. {\it Tail asymptotics for the supremum of a random
walk when the mean is not finite}. Queueing Systems, 46 (2004), 15-33.
\end{thebibliography}
\end{document}